\documentclass[a4paper,11pt]{amsart}

\usepackage{amsmath,amsthm,amsfonts,amssymb,amscd}

\input{xy}

\xyoption{all}

\newtheorem{theorem}{Theorem}[section]
\newtheorem{proposition}[theorem]{Proposition}

\newtheorem{lemma}[theorem]{Lemma}
\newtheorem{corollary}[theorem]{Corollary}
\newtheorem{definition}[theorem]{Definition}
\newtheorem{remark}[theorem]{Remark}

\newenvironment{sis}{\left\{\begin{aligned}}{\end{aligned}\right.}

\numberwithin{equation}{section}

\renewcommand{\P}{\mathbb{P}}
\renewcommand{\O}{\mathcal{O}}
\newcommand{\C}{\mathcal{C}}
\newcommand{\F}{\mathcal{F}}
\newcommand{\M}{\mathcal{M}}
\renewcommand{\L}{\mathcal{L}}
\newcommand{\tD}{\widetilde{\mathcal{D}}_{2g+2}}
\newcommand{\tH}{\widetilde{\mathcal{H}}_g}

\newcommand{\GG}{\mathcal{G}}
\newcommand{\T}{\mathcal{T}_{a,b}}
\newcommand{\tT}{\widetilde{\mathcal{T}}_{a,b}}
\newcommand{\Z}{\mathbb{Z}}

\newcommand{\G}{\mathbb{G}_m}

\newcommand{\Pic}{{\rm Pic}}

\newcommand{\Cl}{{\rm Cl}}

\renewcommand{\S}{{\rm Sym}^{2g+2}(\P^1)}

\newcommand{\B}{\mathbb{B}(2,2g+2)}

\newcommand{\Bsm}{\mathbb{B}_{sm}(2,2g+2)}

\newcommand{\A}{\mathbb{A}(2,2g+2)}

\newcommand{\Asm}{\mathbb{A}_{sm}(2,2g+2)}

\newcommand{\BS}{\mathbb{B}(2,6)}

\newcommand{\BSsm}{\mathbb{B}_{sm}(2,6)}

\renewcommand{\H}{\mathcal{H}_g}

\newcommand{\D}{\mathcal{D}_{2g+2}}

\renewcommand{\det}{{\rm det}}

\begin{document}

\title{Picard group of moduli of hyperelliptic curves}

\author{Sergey Gorchinskiy and Filippo Viviani}

\address{Steklov Mathematical Institute\\ Gubkina str. 8\\ 119991\\ Moscow\\}

\email{serge.gorchinsky@rambler.ru}

\address{Dipartimento di Matematica\\ Universit\`a di Roma Tor Vergata\\
 via della Ricerca Scientifica 1\\ I-00133\\ Roma\\ Italy}

\email{viviani@axp.mat.uniroma2.it}


\thanks{The first author was supported by RFFI grants no. 04-01-00613, 05-01-00455
and INTAS grant no. 05-1000008-8118.}

\thanks{The second author was supported by a grant from the Mittag-Leffler
Institute of Stockholm.}

\maketitle

\begin{abstract}
The main subject of this work is the difference between the coarse moduli 
space and the stack of hyperelliptic curves. We compute their Picard groups, 
giving explicit description of the generators. We get an application to the 
non-existence of a tautological family over the coarse moduli space.
\end{abstract}

\section{Introduction}

Throughout this paper we work over a field $k$ of characteristic different from $2$
and we fix an integer $g\ge 2$. A hyperelliptic curve of genus
$g$ over an algebraically closed field is a smooth curve of genus $g$ which
is a double cover of the projective line $\P^1$ ramified at $2g+2$ points. 
We say that a smooth morphism $f:\F\to S$ of $k$-schemes is a family of 
hyperelliptic curves or of $\P^1$ if any geometric fiber of $f$ is isomorphic 
to a hyperelliptic curve or to $\P^1$, respectively.

In this article we are interested in comparing the (coarse) moduli space $H_g$
of hyperelliptic curves and the moduli stack $\mathcal{H}_g$ of hyperelliptic
curves. 
The stack $\mathcal{H}_g$
has been studied by Arsie and Vistoli (see \cite{AV}
and also \cite{Vis} for $g=2$) who provided a description of it as a quotient stack
and computed its Picard group, which turns out to be isomorphic to
$\Z/(4g+2)\Z$ for $g$ even, and to $\Z/2(4g+2)\Z$ for $g$ odd.
After some auxiliary results on $H_g$ in section \ref{section1}, 
we compute in section \ref{section2} (away from
some bad characteristics of the base field $k$) the class group $\Cl(H_g)$ and 
compare it with $\Pic(\H)$ in Theorem \ref{D-pic} and Corollary \ref{picmoduli}.
As an application we prove in Theorem \ref{non-existence}
the non-existence of a tautological family over $H_g^0$ for $g$ odd, where
$H_g^0$ is the locus of hyperelliptic curves without extra-automorphisms
(for $g$
even a tautological family does not exist over any Zariski open subset in $H_g$).
Also we obtain
that the Picard group of the normal variety $H_g$ is trivial, see
Corollary \ref{picmoduli2}.

Further, for $g=2$ Vistoli proved in \cite{AV} that the Picard group
$\Pic(\mathcal{H}_2)$ is generated by the first Chern class of the Hodge bundle.
In Theorem \ref{gen} from section \ref{section3} we provide an explicit
functorial description of a generator
of the Picard group of the stack $\H$ for arbitrary $g$. Moreover, in Theorem
$\ref{elements}$ we consider some natural elements of the Picard group (obtained
by pushing-forward linear combinations of the relative canonical divisor and the
Weierstrass divisor and then taking the determinant) and express them in terms
of the generator found above. In particular, we show in Corollary \ref{Hodge}
that the first Chern class of the Hodge bundle generates the Picard group if and
only if $4$ does not divide $g$. Otherwise the Hodge line bundle generates
a subgroup of index $2$.

Quite recently Cornalba has computed in \cite{Cor} the Picard
group of
the stack of stable hyperelliptic curves. The article \cite{Cor} also contains a very
beautiful proof of the first assertion of Theorem \ref{gen} over $\mathbb{C}$ by
a quite different method from the one used in the present paper.

Let us finally mention that a more detailed version of this text can be found
at the web in \cite{GV}.

We are grateful to professor A. Ragusa, professor O. Debarre, and professor 
L. Caporaso for organizing an
excellent summer school ``Pragmatic-2004'' held at the University of Catania,
where the two authors began their joint work on this subject.
We thank professor L.~Caporaso who suggested during this summer school
an interesting research problem, from which this work was originated.
We are also grateful to the referee for many useful comments and remarks.

\section{Auxiliary results on the moduli space of 
hyperelliptic curves}\label{section1}

Recall that the coarse moduli space $H_g$ parameterizing isomorphism classes of
hyperelliptic curves is an irreducible variety of dimension $2g-1$ and can be
realized as follows:
\begin{equation}\label{descI}
H_g=(\S-\Delta)/PGL_2,
\end{equation}
where the action of $PGL_2$ is induced from the natural action on $\P^1$
and $\Delta$ is the closed subset in $\S$ consisting of all $(2g+2)$-tuples
on $\P^1$ with at least one coincidence.
We identify $\S$ with the projective space $\B$ of degree $2g+2$ binary forms.
Under this identification $\S-\Delta$ corresponds to the open subset $\Bsm$ of
smooth forms (i.e., whose all roots are distinct) and the action of $PGL_2$ is
defined by the formula $[A]\cdot [f(x)]= [f(A^{-1}x)]$, where $[A]$ is the
class in $PGL_2$ of a $(2\times 2)$ non-degenerate matrix $A$.
By $H_g^0$ denote the open subset
of $H_g$ corresponding to hyperelliptic curves without extra-automorphisms apart
from the hyperelliptic involution. Let $\Bsm^0$ denote the preimage of $H_g^0$
in $\Bsm$. Equivalently, $\Bsm^0$ consists of points in $\Bsm$ with
trivial stabilizers in $PGL_2$.

\begin{proposition}\label{dimension}
The locus $H_g^{aut}=H_g- H_g^0$ has dimension $g$ and hence
codimension $g-1$ in $H_g$. Moreover, it has a unique irreducible component of maximal
dimension corresponding to hyperelliptic curves that have an extra-involution
(besides the hyperelliptic one),
acting on $2g+2$ ramification points as a product of
$g+1$ commuting transpositions.
\end{proposition}

\begin{proof}
The automorphism group ${\rm Aut}(C)$ of a hyperelliptic curve $C$ always
contains the hyperelliptic involution $i$ as a central element. Consider the
group $G={\rm Aut}(C)/\langle i \rangle$. There is a canonical inclusion inside
the symmetric group $G\subset S_{2g+2}$, since every automorphism of a
hyperelliptic curve acts on the ramification divisor. Hence the variety
$H_g^{aut}$ decomposes into the strata
$$
H_g^{aut}=\bigcup_{p\le 2g+2}H_g^{aut,p},
$$
where the union is taken over all primes $p$, $p\le 2g-2$ and
$H_g^{aut,p}$ is the set of hyperelliptic curves such
that there exists an element of order $p$ in the corresponding group
$G$. There is a canonical finite map $H_g^{aut,p-fixed}\to
H_g^{aut,p}$, where $H_g^{aut,p-fixed}$ is the coarse moduli space of
pairs $(C,\sigma)$ such that $C$ is a curve from $H_g^{aut,p}$ and $\sigma$
is an element of order $p$ in the group $G$ associated with $C$.

Since $\sigma\in G$ is uniquely determined by any automorphism of $\P^1$
preserving the ramification divisor, we see that, in fact,
$H_g^{aut,p-fixed}$ is the coarse moduli space of
pairs $(D,\tau)$ such that $D$ is a reduced effective divisor
of degree $2g+2$ on $\P^1$ and $\tau$ is an automorphism of $\P^1$ of order $p$
that satisfies $\tau(D)=D$.

Consider the natural quotient map
$$
\pi:\P^1=\P^1_1\stackrel{p:1}\longrightarrow\P^1_2=\P^1/\langle
\tau\rangle.
$$
Since $p$ is prime, it is well-known that $\pi$ has a
cyclic ramification of order $p$ at two points $x_1,x_2\in\P^1_1$ and $\tau$
is uniquely determined by the points $x_1$ and $x_2$.
There are three possibilities for the divisor
$D\subset \P^1_1$:
\begin{itemize}
\item[0)] $D$ contains no points among $x_1$ and $x_2$,
\item[1)] $D$ contains only one point among $x_1$ and $x_2$,
\item[2)] $D$ contains both points $x_1$ and $x_2$.
\end{itemize}

Hence we get one more stratification:
$$
H_g^{aut,p-fixed}=\bigcup_{l=0,1,2}{H_g^{aut,p-fixed,l}},
$$
where $H_g^{aut,p-fixed,l}$ is
the coarse moduli space of pairs $(R,E)$ such that $R$ and $E$ are
non-intersecting reduced effective divisors on $\P^1_2$ of degrees $2$ and
$(2g+2-l)/p$, respectively (in particular, we require that $2g+2-l$
is divisible by $p$). Thus we get the equality
$$
\dim H_g^{aut,p-fixed,l}=2+\frac{2g+2-l}{p}-3=\frac{2g+2-l}{p}-1.
$$
Notice that the case $p=2$ and $l=1$ is impossible because of
the divisibility condition. Further, if $p\ge 3$, then
$$
\frac{2g+2-l}{p}-1\le \frac{2g+2}{3}-1\le g-1
$$
and for $p=2$, $l=2$ we have
$$
\frac{2g+2-2}{2}-1=g-1.
$$
So, we get the inequality
$$
\dim\left(H_g^{aut,2-fixed,2}\cup
\bigcup_{3\le p} H_g^{aut,p-fixed}\right)=
\max_{(p,l)\ne(2,0)} \{\dim (H_g^{aut,p-fixed,l})\}\le g-1.
$$

Suppose that $p=2$ and $l=0$; then $\dim(H_g^{aut,2-fixed,0})=g$.
We claim that in this case the curve $C$ has an
element $\tilde{\sigma}$ of order two in the automorphism group ${\rm Aut}(C)$ itself
(not only in $G$). Indeed, consider the composition
$$
\varphi:C\stackrel{2:1}\longrightarrow
\P^1_1\stackrel{2:1}\longrightarrow\P^1_2.
$$
This is a Galois map of degree $4$ with Galois group $H$ generated
in ${\rm Aut}(C)$ by any preimage $\tilde{\sigma}\in {\rm Aut}(C)$ of $\sigma\in
G$ and $i$. It is easily seen that the ramification of
$\varphi$ consists only of pairs of double points. If $H\cong
\mathbb{Z}/4\mathbb{Z}$, then the inertia group at all
ramification points of $\varphi$ would be the same, namely $\langle i \rangle\subset H$.
This would mean that the map $\pi\colon \P^1_1=C/\langle i
\rangle\to\P^1_2$ is unramified. This contradiction shows that
$H\cong \mathbb{Z}/2\mathbb{Z}\times \mathbb{Z}/2\mathbb{Z}$
and $\tilde{\sigma}\in {\rm Aut}(C)$ has order two.

Conversely, if ${\rm Aut}(C)$ has an element $\sigma\ne i$ of order two,
then the corresponding number $l$ equals zero. Indeed, otherwise the inertia
group of $\varphi$ at any point from
$D\cap\{x_1,x_2\}$ would be isomorphic to $\Z/4\Z$, hence $H$ would be
isomorphic to $\Z/4\Z$ and would have only one element
of order two.

Note that $H_g^{aut,2-fixed,0}$ is irreducible and, moreover, it follows from
the explicit geometric description of the ramification of the
covering $\varphi:C\to \P^1_2$ that $\sigma\in G\subset
S_{2g+2}$ must be equal to the product of $g+1$ commuting transpositions.
This completes the proof.
\end{proof}

\begin{remark}
It is possible to give a purely combinatorial proof of a weaker version of
this proposition (see \cite[Prop. 4.3']{GV}).
\end{remark}

\begin{remark}
It follows from Proposition \ref{dimension}
that for $g\geq 3$ the smooth locus of the normal variety $H_g$ is equal to $H_g^0$
(see \cite[Proposition 4.5]{GV}).
\end{remark}

\begin{remark}
It is interesting to compare the above results with the analogous ones for
the coarse moduli space $M_g$ of smooth curves of genus $g\geq 3$.
The locus $M_g^{aut}$ of curves with non-trivial automorphisms
is a closed subset of dimension $2g-1$, and it has a unique irreducible
component of maximal dimension corresponding to hyperelliptic curves.
Moreover, the smooth locus of $M_g$ is equal to $M_g^0$ if $g\geq 4$,
while the smooth locus of $M_3^{smooth}$ is equal to $M_3^0\cup H_3^0$
(see \cite{Rau}, \cite{Pop},\cite{Oort}, \cite{Lon2}).
\end{remark}

The following result is needed for the sequel.

\begin{lemma}\label{g=2}
Let $D$ be the unique irreducible divisor on $\BSsm$ from
$\BSsm-\BSsm^0$ (see Proposition \ref{dimension}) and let $\overline{D}$
be its closure in $\BS$. Then $\overline{D}$ is an irreducible
hypersurface in $\BS=\P^6$ of degree $15$.
\end{lemma}

\begin{proof}
Consider the natural map
$\pi:(\P^1)^6\stackrel{/S_6}{\longrightarrow}
{\rm Sym}^6(\P^1)$. Suppose that an element $\sigma\in S_6$ is conjugate
to the permutation $(12)(34)(56)$. Denote by $D_{\sigma}$ the divisor in $(\P^1)^6$
consisting of all points $(P_1,\ldots,P_6)\in (\P^1)^6$ such that there exists
a non-trivial element $A\in PGL_2$ that satisfies
$A(P_1,\ldots,P_6)=(\sigma(P_1),\ldots,\sigma(P_6))$.

It follows from Proposition \ref{dimension} that
$\pi^{-1}(\overline{D})=\bigcup_{\sigma}
D_{\sigma}$, where the union is taken over the $15$
elements of $S_6$ conjugated to $(12)(34)(56)$.

Let us compute the class of $D_{\sigma}$ in the
Picard group $\Pic((\P^1)^6)\cong(\mathbb{Z})^6$. Without loss of
generality we may suppose that $\sigma=(12)(34)(56)$. Take a line
$l=\{P_1\}\times\ldots\times\{P_5\}\times\P^1$ in $(\P^1)^6$ for
general points $P_i\in\P^1$. It is well-known that there exists
a unique non-trivial element $A\in PGL_2$ exchanging $P_1$ with
$P_2$ and $P_3$ with $P_4$. In particular, $A$ has order two. Hence
the point $P_6=A(P_5)$ is uniquely determined and the intersection
$l\cap D_{(12)(34)(56)}$ consists of
one point. It is easy to prove that this intersection is actually
transversal. By the symmetry of $D_{(12)(34)(56)}$, the same is
true for all other ``coordinate'' lines in $(\P^1)^6$ and the
class of $\overline{D}_{(12)(34)(56)}$ in $\Pic((\P^1)^6)$ is equal
to $(1,1,1,1,1,1)$. This completes the proof of Lemma \ref{g=2}.
\end{proof}

\section{Comparison between Picard groups of moduli space and stack of
hyperelliptic curves}\label{section2}

Recall that $\H$ is a category such that objects in $\H$ are families
$\pi:\F\rightarrow S$ of
hyperelliptic curves of genus $g$ and morphisms in $\H$ are Cartesian diagrams
between such families. Associating the base to a family, we obtain that $\H$ is a
category fibered in groupoids over the category of $k$-schemes.
By $\H^0$ denote the full fibered subcategory of $\H$ such that the
objects in $\H^0$ are families of hyperelliptic curves whose geometric fibers
have no extra-automorphisms.

Let us cite from \cite{AV} two fundamental facts about the fibered category $\H$. We
keep notations from the previous section.

\begin{theorem}[Arsie--Vistoli] \label{stack}

\hspace{0cm}

\begin{enumerate}

\item[(i)]
The fibered category $\H$ is a Deligne--Mumford
algebraic stack and can be realized as $\mathcal{H}_g=[\Asm/(GL_2/\mu_{g+1})]$
with the usual action given by $[A]\cdot f(x)=f(A^{-1}\cdot x)$, where $\mu_n$
denotes the group scheme of $n$-th power roots of unity for $n$ prime to ${\rm char}(k)$.

\item[(ii)]
Suppose that ${\rm char}(k)$ does not divide $2g+2$.
Then the Picard group $\Pic(\mathcal{H}_g)$ is the 
quotient of $\Pic^{GL_2/\mu_{g+1}}(\A)=(GL_2/\mu_{g+1})^*=\Z$ 
of order equal to $4g+2$ if $g$ is even and $2(4g+2)$ if $g$ is odd
(where $G^*={\rm Hom}(G,\G)$ for an algebraic group $G$).

\end{enumerate}

\end{theorem}
In addition, there is a well-known explicit description of
$(GL_2/\mu_{2g+2})^*$.

\begin{lemma}\label{quotientgroup}

\hspace{0cm}

\begin{itemize}

\item[(i)] If $g$ is even, then there is an isomorphism of algebraic groups
$GL_2/\mu_{g+1}$ $\to GL_2$ given by $[A]\mapsto {\rm det}(A)^{\frac{g}{2}}A$
and $(GL_2/\mu_{g+1})^*=\Z$ is generated by ${\rm det}^{g+1}$.

\item[(ii)] If $g$ is odd, then there is an isomorphism of algebraic groups
$GL_2/\mu_{g+1}$ $\to \mathbb{G}_m\times PGL_2$ given by
$[A]\mapsto ({\rm det}(A)^{\frac{g+1}{2}},[A])$ and
$(GL_2/\mu_{g+1})^*=\Z$ is generated by ${\rm det}^{\frac{g+1}{2}}$.

\end{itemize}

\end{lemma}

Now we compare the stack $\H$ and its coarse moduli scheme $H_g$. In particular,
we compare the Picard group of $\H^0$ and the Picard group of $H_g^0$. With this
aim it is natural to introduce a new stack, which is ``intermediate''
between $\H$ and $H_g$.

\begin{definition}\label{newfunctor}
Let $\mathcal{D}_{2g+2}$ be a category such that the objects in $\D$
are families $p:\C\to S$ of $\P^1$ together with an effective Cartier divisor $D\subset \C$
finite and \'etale over $S$ of degree $2g+2$ and and the morphisms in $\D$
are natural Cartesian diagrams.
\end{definition}

Associating the base to each family, we obtain that $\D$ is
a category fibered in groupoids over the category of $k$-schemes. We say that a
divisor $D\subset \P^1$ has no automorphisms if there is no non-trivial element
$f\in{\rm Aut}(\P^1)$ such that $f(D)=D$. By $\D^0$ denote the full fibered
subcategory of $\D$ such that the objects in $\D^0$ are families whose all
geometric fibers have no automorphisms.

The following result is analogous to theorem \ref{stack}, as well as its proof.

\begin{proposition}\label{D-stack}
The fibered category
$\D$ is a Deligne--Mumford algebraic stack and can be realized as $[\Bsm/PGL_2]$
with the usual action given by
$[A]\cdot [f(x)]=[f(A^{-1}\cdot x)]$. Moreover, there is a natural
isomorphism of stacks $[\Bsm/PGL_2]\cong [\Asm/(GL_2/\mu_{2g+2})]$, where we consider
the same action of $GL_2$ on $\Asm$ as in theorem \ref{stack}.
\end{proposition}

\begin{proof}
Consider the auxiliary functor $\tD$ that associates
with a $k$-scheme $S$ the set of collection
$$
\tD(S)=\{(\C\stackrel{p}{\to} S, D, \phi:\C\cong \P^1_S)\},
$$
where the family $\C\to S$ and the divisor $D$ are as in Definition \ref{newfunctor}
and $\phi$ is an isomorphism over $S$ between the family $\C\to S$ and the trivial family
$\P^1_S:=S\times_k \P^1$.
Clearly, $\tD\cong{\rm Hom}(-,\Bsm)$. The group sheaf $\underline{\rm Aut}(\P^1)=
PGL_2$ acts on $\tD$ by composing with the isomorphism
$\phi$ and it is easy to check that the corresponding action of $PGL_2$ on
$\Bsm$ is the one given in the statement of Proposition \ref{D-stack}.
Finally, descent theory implies that the
forgetful morphism $\tD\to \D$ is a principle bundle over $\D$
with the group $PGL_2$. Thus we get
the description of $\D$ as a quotient stack $[\Bsm/PGL_2]$.

To prove the second part of the proposition observe that, applying Lemma
\ref{quotientgroup}(ii) with $g+1$ replaced by $2g+2$, one deduces an
isomorphism $GL_2/\mu_{2g+2}\cong \G\times PGL_2$. It is easily shown
that the corresponding action of
$\G\times PGL_2$ on $\Asm$ is given by
$(\alpha, [A])\cdot f(x)= \alpha^{-1}\cdot ({\rm det}A)^{g+1}f(A^{-1}\cdot x)$.
Hence the quotient stack of $\Asm$ by $GL_2/\mu_{2g+2}\cong \G\times PGL_2$
can be taken in two steps: first, we take the quotient over the subgroup
$\G/\mu_{2g+2}\cong\G$,
which is isomorphic to $\Bsm$ since the action is free, and then we take the
quotient over
$GL_2/\G\cong PGL_2$ with the usual action.
\end{proof}

From these explicit descriptions we get a diagram
$$
\xymatrix{
\mathcal{H}_g \ar[rr]^{\Psi} \ar[rd]_{\Phi_{\mathcal{H}}} & &
\mathcal{D}_{2g+2}
\ar[ld]^{\Phi_{\mathcal{D}}}\\
&  H_g & }\\
$$
where $H_g$ is the coarse moduli space for both stacks and the morphism
$\Psi:\H\to \D$ corresponds to the fact that every family  $\pi:\F\to S$ of
hyperelliptic curves is
a double cover of a family $p:\C\to S$ of $\P^1$ such that the ramification divisor
$W\subset\F$
and the branch divisor $D\subset \C$ are both finite and \'etale over $S$ of degree
$2g+2$ (see \cite{LK}).

The following result is well-known.

\begin{corollary}\label{universaldivisor}
There is an isomorphism $\D^0\cong H_g^0$, i.e.,
$H_g^0$ is a fine moduli scheme for $\mathcal{D}_{2g+2}^0$.
\end{corollary}

\begin{proof}
By Proposition \ref{D-stack}, $\D^0\cong[\Bsm^0/PGL_2]$. By definition,
the action of $PGL_2$ on $\Bsm^0$ is free and a standard fact is that that the
line bundle $\O(1)$ on $\B$ admits a canonical $PGL_2$-linearization.
Therefore, $[\Bsm^0/PGL_2]=\Bsm^0/PGL_2=H_g^0$.
\end{proof}

Now we compute the Picard group of the stack $\D$ and relate it to
the Picard group of $\H$.

\begin{theorem}\label{D-pic}
The natural map $\Pic(\D)\to \Pic(\H)$ is injective, being
an isomorphism for $g$ even and an inclusion of index $2$ for $g$ odd.
Therefore if ${\rm char}(k)$ does not divide $2g+2$ then $\Pic(\D)=\Z/(4g+2)\Z$.
\end{theorem}

\begin{proof}
Theorem \ref{stack}(i) and the second description from Proposition
\ref{D-stack} show that there exists a Cartesian diagram
$$
\begin{array}{ccc}
\H&\to&[\A/(GL_2/\mu_{g+1})]\\
\downarrow&&\downarrow\\
\D&\to&[\A/(GL_2/\mu_{2g+2})],
\end{array}
$$
where each horizontal arrow is an open embedding.
Hence the statement in question reduces to the study of the map
$\Pic^{GL_2/\mu_{2g+2}}(\A)=(GL_2/\mu_{2g+2})^*\to \Pic^{GL_2/\mu_{g+1}}(\A)=
(GL_2/\mu_{g+1})^*$. The needed results follow immediately from Lemma
\ref{quotientgroup} and Theorem \ref{stack}(ii).
\end{proof}

\begin{remark}\label{comparison}
It can easily be checked that the map
$\Pic^{PGL_2}(\B)=$$\langle\O(1)\rangle$$  \to
\Pic^{GL_2/\mu_{g+1}}(\A)$ takes $\O(1)$ to the character
defined by $\det^{g+1}$.
\end{remark}

\begin{corollary}\label{picmoduli}
Suppose that ${\rm char}(k)$ does not divide $2g+2$ and neither is equal to $5$ if $g=2$.
Then we have
$$
\Cl(H_g)=\Cl(H_g^0)=\Pic(H_g^0)=\Pic(\D^0)=\left\{
\begin{aligned}
&\mathbb{Z}/(4g+2)\mathbb{Z} \text{ if } g \geq 3, \\
&\mathbb{Z}/5\mathbb{Z} \text{ if } g=2.\\
\end{aligned}
\right.
$$
Moreover the natural map $\Pic(H_g^0)\to \Pic(\mathcal{H}_g^0)$ is injective, being
an isomorphism for $g$ even and an inclusion of index $2$ for $g$ odd. 
\end{corollary}

\begin{proof}
By Theorem \ref{stack}, $\mathcal{H}_g^0\cong[\Asm^0/(GL_2/\mu_{g+1})]$.
Since $H_g^0$ is smooth, $\Cl(H_g^0)=\Pic(H_g^0)$. By Corollary \ref{universaldivisor},
$\Pic(H_g^0)=\Pic(\D^0)$. By Proposition \ref{dimension},
we have $\Pic(\mathcal{H}_g)=\Pic(\mathcal{H}_g^0)$,
$\Pic(\mathcal{D}_{2g+2})=\Pic(\mathcal{D}_{2g+2}^0)$, and $\Cl(H_g)=\Cl(H_g^0)$
when $g-1\ge 2$. Thus the needed
statement for $g\ge 3$ follows from Theorem \ref{stack}(ii) and Theorem \ref{D-pic}.
For $g=2$ the Cartesian diagram
$$
\begin{array}{ccc}
\mathcal{H}_2^0&\to&\mathcal{H}_2\\
\downarrow&&\downarrow\\
\mathcal{D}_2^0&\to&\mathcal{D}_2\\
\end{array}
$$
shows that the natural map $\Pic(\mathcal{D}^0_2)\to\Pic(\mathcal{H}_2^0)$
is an isomorphism. Further, it follows from Lemma \ref{g=2}, Theorem \ref{stack}(ii),
and Remark \ref{comparison} that $\Pic(\mathcal{H}_2^0)=\Z/5\Z$.

On the other hand, by the result of Igusa (see \cite{Igu}), for ${\rm char}(k)\neq 5$
we have $H_2=\mathbb{A}^3/\mu_5$, where $\zeta\in\mu_5$ acts on $\mathbb{A}^3$ by
formula $(x_1,x_2,x_3)\mapsto (\zeta x_1,\zeta^2 x_2,\zeta^3 x_3)$.
Being the image of the origin,
the unique singular point of $H_2$ corresponds to the curve
$C_0:=\{y^2=x^6-x\}$. Therefore,
$$
\Cl(H_2)=\Cl(H_2-[C_0])=\Pic(H_2-[C_0])=
\Pic^{\Z/5\Z}(\mathbb{A}_k^3-\{0\})\cong \Z/5\Z,
$$
and the natural surjective morphism $\Cl(H_2)\to\Cl(H_2^0)=\Pic(H_2^0)=
\Pic(\mathcal{H}_2^0)$ is actually an isomorphism.
\end{proof}

\begin{remark}
It is possible to compute the Picard group $\Pic(H_g^0)$ directly, without using the stack
description from Corollary \ref{universaldivisor} and the results of Theorem
\ref{stack} (see \cite[Theorem 4.7]{GV}).
\end{remark}

\begin{corollary}\label{picmoduli2}
Suppose that ${\rm char}(k)$ does not divide $(2g+1)(2g+2)$.
Then\\ $\Pic(H_g)=0$.
\end{corollary}

\begin{proof}
Since $H_g$ is a normal variety, the map $\Pic(H_g)\to \Cl(H_g)$ is
injective. Put $N_2=5$ and $N_g=4g+2$ if $g\ge 3$.
By Corollary \ref{picmoduli},
$\Cl(H_g)=\Pic(\D^0)$ is a cyclic group of order
$N_g$ generated by the image of the character
$\det^{g+1}$ under the natural
map $\Pic^{GL_2/\mu_{2g+2}}(\A)=(GL_2/\mu_{2g+2})^*\to
\Pic^{GL_2/\mu_{2g+2}}(\Asm^0)=\Pic(\D^0)$.
Hence $\Pic(H_g)$ is contained inside the subgroup of $\Pic(\D)$
generated by the images of characters $\chi\in(GL_2/\mu_{2g+2})^*$ such that
the restriction of $\chi$ to the $GL_2$-stabilizer of any point in $\Asm$
is equal to a multiple of the character $\det^{N_g(g+1)}$.

Since ${\rm char}(k)$ does not divide $(2g+1)(2g+2)$, the binary forms
$f_1:=X^{2g+1}Y-Y^{2g+2}$ and $f_2:=X^{2g+2}-Y^{2g+2}$ belong to $\Asm$. The
$GL_2$-subgroups ${\rm diag}\{\mu_{2g+1},1\}$ and
${\rm diag}\{\mu_{2g+2},1\}$ stabilize $f_1$ and $f_2$, respectively.
This concludes the proof.
\end{proof}

\begin{remark}
For the moduli spaces of smooth curves of genus $g\geq 3$ over $\mathbb{C}$
we have
$$
\Z=\Pic(M_g)\hookrightarrow
\Cl(M_g)=\Cl(M_g^0)=\Pic(M_g^0)=\Pic(\mathcal{M}_g^0)=\Pic(\mathcal{M}_g)=\Z.
$$
Still the index of the first
group inside the second one remains unknown (see \cite[section 4]{AC}).
\end{remark}

Now let us give an application of the comparison between Picard groups.
Recall that a tautological family of hyperelliptic curves exists over a
non-empty Zariski open subset in $H_g$ if and only if $g$ is odd
(see \cite[Exercise 2.3]{HM},
where ``universal'' should be replaced by ``tautological'').
For $g$ odd we get the following non-existence result.

\begin{theorem}\label{non-existence}
For $g$ odd there does not exist a tautological family over $H_g^0$ (and
henceforth over all $H_g$).
\end{theorem}

\begin{proof}
A tautological family over $H_g^0$ would define a section of the modular map
$H_g^0\to \mathcal{H}_g^0$.
Consequently there would be a
splitting of the map $\Pic(H_g^0)\hookrightarrow \Pic(\mathcal{H}_g^0)$, which
is impossible because of the explicit description of this map in Corollary
\ref{picmoduli}.
\end{proof}

\begin{remark}
It is possible to give a direct proof of the last statement: first, compute directly the
Picard group of the universal family of $\P^1$ over $H_g^0$, and then find
explicitly the class of the
universal divisor in this group and check that it is not divisible by two if
$g$ is odd, see \cite[Proposition 6.13]{GV}.
\end{remark}

\begin{remark}
For $g=1$ there exists a tautological family
over $H_1^0$ (in this case $\Pic(H_1^0)=0$), see \cite[page 58]{Mum}.
\end{remark}

\section{Explicit generators of the Picard group}\label{section3}

In this section we give an explicit construction for a generator of the Picard group
${\rm Pic}(\H)$ in terms of Mumford's functorial description of the Picard group
of a stack (see \cite{Mum}, \cite{EG}).

Let $\pi:\F\to S$ be a family of hyperelliptic curves. In the discussion after
Proposition \ref{D-stack} we introduced a family $p:\C\to S$ of $\P^1$ and two
Cartier divisors $W\subset \F$ (so called Weierstrass divisor) and $D\subset
\C$. By classical theory of double covers, there exists a line bundle
$\L$ on $\C$ such that
$(\L^{-1})^{\otimes 2}=\O_{\C}(D)$. This line bundle satisfies two relations:
\begin{equation}\label{f2}
f^*(\L^{-1})=\O_{\F}(W),
\end{equation}
\begin{equation}\label{f3}
f_*(\O_{\F})=\O_{\C}\oplus \L.
\end{equation}
Moreover, Hurwitz formula tells that
\begin{equation}\label{f4}
\omega_{\F/S}=f^*(\omega_{\C/S})\otimes \O_{\F}(W).
\end{equation}
One can check that the fibered category $\H$ is equivalent to the fibered category $\H'$
such that the objects in $\H'$ are collections
$(\C\stackrel{p}{\to} S, \L, \L^{\otimes 2}\stackrel{i}{\hookrightarrow}
\O_{\C})$, where $p:\C\to S$ is a family of $\P^1$, $\L\in\Pic(\C)$ and the morphisms in 
$\H'$ are natural Cartesian diagrams (see \cite[section 2]{AV}).

\begin{theorem}\label{gen}
Let $\GG$ be the element in $\Pic(\H)$ such that
for any family of hyperelliptic
curves $\pi:\F\to S$ with the Weierstrass divisor $W$ the line bundle
$\GG(\pi)$ on $S$ is defined by the formula
$$
\GG(\pi)=\begin{sis}
\pi_*\left(\omega_{\F/S}^{-(g+1)}\left((g-1)W\right)\right)& \text{ if } g \text {
is even, }\\
\pi_*\left(\omega_{\F/S}^{-\frac{(g+1)}{2}}\left(\frac{g-1}{2}W\right)\right) &
\text { if } g \text { is odd. }\\
\end{sis}
$$
Then $\GG$ generates the Picard group $\Pic(\H)$ and equals to the image of the
character $\chi_0$ under the natural map
$\Pic^{GL_2/\mu_{g+1}}(\A)=(GL_2/\mu_{g+1})^*$ $ \to\Pic(\H)$,
where $\chi_0=\det^{g+1}$ for $g$
even and $\chi_0=\det^{\frac{g+1}{2}}$ for $g$ odd.
\end{theorem}

\begin{proof}
By Theorem \ref{stack} and Lemma \ref{quotientgroup},
${\rm Pic}(\H)={\rm Pic}^{GL_2/\mu_{g+1}}(\Asm)$ is a
cyclic group generated by the trivial line bundle $\Asm\times k$, on which
$GL_2/\mu_{g+1}$ acts by the character $\chi_0$.

Following the proof of Theorem \ref{stack} from \cite{AV},
consider the auxiliary functor $\tH$ that associates with a $k$-schemes $S$ the
set of collections
$$
\tH(S)=\{(\C\stackrel{p}{\to} S, \L, \L^{\otimes
2}\stackrel{i}{\hookrightarrow} \O_{\C}, \phi:(\C,\L)\cong (\P^1_S,\O_{\P^1_S}
(-g-1)))\},
$$
where $p:\C\to S$, $\L\in\Pic(\C)$ are as
above and the isomorphism $\phi$ consists of an isomorphisms of $S$-schemes
$\phi_0:\C\cong \P^1_S$ plus an isomorphism of invertible sheaves
$\phi_1:\L\cong \phi_0^*\O_{\P^1_S}(-g-1)$.

In \cite{AV}, it was proved that $\tH\cong {\rm Hom}(-,\Asm)$
and that the forgetful morphism $\tH\to \H'\cong \H$ is a principal bundle
over $\H$ with the group $GL_2/\mu_{g+1}$.

Consider the following commutative diagram of $GL_2/\mu_{g+1}$-equivariant maps:
$$
\xymatrix{
\tH\times k \ar[r]^(0.35){\cong}\ar[d] & \Asm\times k\ar[d]\\
\tH\ar[r]^(0,35){\cong} & \Asm.
}
$$

The functor $\tH\times k$ associates with a $k$-scheme $S$ the set of
collections
$$
(\tH\times k)(S)=\left\{(\C\stackrel{p}{\to} S, \L,
\L^{\otimes 2}
\stackrel{i}{\hookrightarrow} \O_{\C}
, \phi:(\C,\L)\cong (\P^1_S,\O_{\P^1_S}(-g-1)),\M)\right\}
$$
where $\M=\O_S$ is the trivial line bundle with the action of the group
$(GL_2/\mu_{g+1})(S)$ given by the character $\chi_0$.

Put $\P^1_S=\P(V_S)$, where $V$ is a two-dimensional vector space over the
ground field $k$. From the Euler exact sequence for the
trivial family $p_S:\P^1_S\to S$
$$
0\to \O_{\P^1_S}\to p_S^*(V_S^*) (1)\to
\omega_{\P^1_S/ S}^{-1}\to 0
$$
one deduces the $(GL_2/\mu_{g+1})(S)$-equivariant isomorphism
\begin{equation}\label{Euler}
p_S^*(({\rm det}V_S)^{-1})\otimes \O_{\P^1_S}(2)\cong
\omega^{-1}_{\P^1_S/S},
\end{equation}
where we consider the canonical actions of $(GL_2/\mu_{g+1})(S)$ on $\P^1_S$ and
on the invertible sheaves involved.
Using projection formula, the fact that $(p_S)_*(\O_{\P^1_S})=\O_S$, and the
$(GL_2/\mu_{g+1})(S)$-equivariant identity $(\det V_S)^{-(g+1)}$ $=\M$ for $g$ even,
and $(\det V_S)^{-\frac{g+1}{2}}=\M$ for $g$ odd, we get
$(GL_2/\mu_{g+1})(S)$-equivariant isomorphisms
$$\begin{sis}
\M&\cong
(p_S)_*\left(\omega_{\P^1_S/S}^{-(g+1)}\otimes \O_{\P^1_S}(-(2g+2))\right)
&\text { if } g \text { is even, }\\
\M&\cong
(p_S)_*\left(\omega_{\P^1_S/S}^{-\frac{g+1}{2}}\otimes \O_{\P^1_S}(-(g+1))\right)
&\text{ \: if } g \text{ is odd. }
\end{sis}$$
Let us remark that $\phi:(\C,\L)\cong (\P^1_S,\O_{\P^1_S}(-g-1))$ induces a
canonical isomorphism $\omega_{\C/S}\cong\omega_{\P^1_S/S}$ defined by the
isomorphism $\phi_0$.
Hence the quotient line bundle
$\GG'=\left[\left(\tH\times
k\right)/\left(GL_2/\mu_{g+1}\right)\right]$
over $\H'$ is isomorphic to
$$
\GG'(p)=
\begin{sis}
&\left\{\C\stackrel{p}{\to} S, \L, \L^{\otimes 2}
\stackrel{i}{\hookrightarrow} \O_{\C}, p_*\left(\omega_{\C/S}^{-(g+1)}\otimes
\L^{2}\right)\right\} &\text { if } g \text { is even, }\\
&\left\{\C\stackrel{p}{\to} S, \L, \L^{\otimes 2}
\stackrel{i}{\hookrightarrow} \O_{\C},
p_*\left(\omega_{\C/S}^{-\frac{g+1}{2}}\otimes
\L\right)\right\} &\text { if } g \text { is odd.}\\
\end{sis}
$$
Now we express the preceding line bundles as push-forwards with respect to the map $f$
of line bundles on the hyperelliptic family $\pi:\F\to S$. Using formulas (\ref{f2}) and
(\ref{f4}), we get
$$
\begin{sis}
&f^*\left(\omega_{\C/S}^{-(g+1)}\otimes \L^{2}\right)=
\omega_{\F/S}^{-(g+1)}\left((g-1)W\right) &\text { if } g \text { is even, }\\
&f^*\left(\omega_{\C/S}^{-\frac{g+1}{2}}\otimes \L\right)=
\omega_{\F/S}^{-\frac{g+1}{2}}\left(\frac{g-1}{2}W\right)& \text{ if } g \text{
is odd. }\\
\end{sis}
$$
To conclude the proof it remains to note that
the line bundles $\omega_{\C/S}^{g+1}\otimes
\L^{\otimes(-2)}$ and $\omega_{\C/S}^{\frac{g+1}{2}}\otimes
\L^{-1}$ for $g$ odd are trivial on each geometric fiber of $p$
and $p_*=\pi_*f^*$ for them.
\end{proof}

Now we express some natural elements of $\Pic(\H)$ in
terms of the generator found above. Recall that
given a family $\pi:\F\to S$ of hyperelliptic curves, there are two natural line
bundles
on $\F$: the relative canonical line bundle $\omega_{\F/S}$ and
the line bundle $\O_{\F}(W)$ associated with the Weierstrass
divisor $W=W_{\F/S}$. Consider their
multiple $\omega_{\F/S}^a\otimes \O_{\F}(bW)$ for any integers $a$ and $b$.
Note that
it restricts to any geometric fiber $F$ of the family $\pi$ as
$$
\omega_{\F/S}^a\otimes \O_{\F}(bW)|_{F}=aK_F+bW_F=[a(g-1)+b(g+1)]g^1_2=
[(a+b)g+b-a]g_2^1.
$$
Since $h^0(F,\O_F(k g^1_2))=k+1$ for any non-negative integer $k$,
the push-forward $\pi_*(\omega_{\F/S}^a\otimes \O_{\F}(bW))$ is a vector
bundle of rank $m(a,b)+1$ on the base $S$ if $m(a,b)\ge 0$, where
$m(a,b):=(a+b)g+b-a$. Let $\T$ be an element in ${\rm Pic}(\H)$ defined by the formula
$$
\T(\pi)={\rm det }\left(\pi_*(\omega_{\F/S}^a\otimes \O_{\F}(bW))\right)\in
{\rm Pic}(S).
$$

\begin{proposition}\label{elements}
If $0\leq m(a,b)< g+1$, then
$$
\T=\begin{sis}
\GG^{-\frac{(a+b)(m(a,b)+1)}{2}}& \text{ if } g \text{ is even, }\\
\GG^{-(a+b)(m(a,b)+1)} & \text{ if } g \text{ is odd, }
\end{sis}
$$
and if $m(a,b)\ge g+1$, then
$$
\T=\begin{sis}
\GG^{\frac{(a+b-1)(g-m(a,b))}{2}}& \text{ if } g \text{ is even, }\\
\GG^{(a+b-1)(g-m(a,b))} & \text{ if } g \text{ is odd. }
\end{sis}
$$
\end{proposition}

\begin{proof}
It follows from formulas (\ref{f2}), (\ref{f3}), and (\ref{f4}) that $\T\in\Pic(\H)$
corresponds to an element $\T'\in\Pic(\H')$ that associates with an object
$(\C\stackrel{p}{\to}S, \L, \L^{\otimes 2}\stackrel{i}{\hookrightarrow}\O_{\C})$ from 
$\H'(S)$ the line bundle
$$
\T'(p)={\rm det} \:p_*\left(\omega_{\C/S}^a\otimes \L^{-(a+b)}\right)\otimes
{\rm det} \:p_*\left(\omega_{\C/S}^a\otimes \L^{-(a+b)+1}\right)\in
{\rm Pic}(S).
$$
Now we compute the pull-back $\tT$ of $\T'$ to ${\rm
Pic}(\tH)$. Using the isomorphism $\phi:(\C,\L)\cong (\P^1_S=\P(V_S),\O(-g-1))$ and the
Euler
formula (\ref{Euler}), we obtain
$$
\omega_{\C/S}^a\otimes \L^{-(a+b)}\cong p_S^*((\det V_S)^a)\otimes \O_{\P^1_S}
(-2a)\otimes
\O_{\P^1_S}\left((a+b)(g+1)\right)=
$$
$$
=p_S^*((\det V_S)^a)\otimes \O_{\P^1_S}(m(a,b))
$$
and, analogously,
$$
\omega_{\C/S}^a\otimes \L^{-(a+b)+1}\cong
p_S^*((\det V_S)^a)\otimes \O_{\P^1_S}(m(a,b)-(g+1)).
$$
Using the relation $\det({\rm Sym}^n(V_S))=(\det
V_S)^{\frac{n(n+1)}{2}}$, we get
$$
\det\:(p_S)_*\left(p_S^*(\det V_S)^a)\otimes \O_{\P^1_S}(m(a,b))\right)=
\det\left((\det V_S)^{a} \otimes {\rm Sym}^{m(a,b)}(V_S)\right)=
$$
$$
={(\det V_S)}^{(a+b)(g+1)\frac{(m(a,b)+1)}{2}},
$$
and, analogously,
$$
\det\:(p_S)_*\left(p_S^*(\det V_S)^a)\otimes \O_{\P^1_S}(m(a,b)-(g+1))\right)=
(\det V_S)^{(a+b-1)(g+1)\frac{m(a,b)-g}{2}}
$$
if $m(a,b)\ge g+1$ and, otherwise, the latter push-forward is zero.
We conclude using the relation $(\det V_S)^{-(g+1)}=\M$, where $\M$ is the pull-back
from $\Pic(\H)$ to $\Pic(\tH)$ of the generator $\GG$ (see the proof
of Theorem \ref{gen}).
\end{proof}

\begin{remark}
It is possible to prove the relations from Proposition \ref{elements} directly
for a given family of hyperelliptic
curves $\pi:\F\to S$ without using this stack computation, see \cite[Theorem 5.8]{GV}.
\end{remark}

Among the elements $\T$ one is of particular interest, namely the Hodge line
bundle $\pi_*(\omega_{\F/S})$, which equals $\mathcal{T}_{1,0}$. The
following result was proved for $g=2$ by Vistoli in \cite{Vis}.
\begin{corollary}\label{Hodge}
\hspace{0cm}
\begin{itemize}
\item[(i)]
The Hodge line bundle is equal to
$$
\det \: \pi_*(\omega_{\F/S})=
\begin{sis}
&\GG^{g/2}& \text{ if } g \text{ is even, }\\
&\GG^{g}&\text{ if } g \text{ is odd. }
\end{sis}
$$
\item[(ii)]
Suppose that ${\rm char}(k)$ does not divide $2g+2$. Then the Hodge line bundle
generates the Picard group $\Pic(\H)$ if $g$ is not divisible by $4$. Otherwise,
it generates a subgroup of index $2$ in $\Pic(\H)$.
\end{itemize}

\end{corollary}

\begin{remark}
The Hodge line bundle generates the Picard group of
$\mathcal{M}_g$ over $\mathbb{C}$ (see \cite{Harer} and \cite{AC}).
\end{remark}

Let us remark that for $g$ even the generator $\GG$ of the Picard group $\Pic(\H)$
equals to $\mathcal{T}_{g/2,1-g/2}={\rm det }\left(\pi_*(\omega_{\F/S}^{g/2}\otimes
\O_{\F}((1-g/2)W))\right)$. It follows from the discussion before
Proposition \ref{elements} that the restriction of the line bundle
$\omega_{\F/S}^{g/2}\otimes \O_{\F}((1-g/2)W)$ to any geometric fiber of the family
$\pi$ is equal to $g^1_2$. Note that, being non-unique,
a line bundle on $\F$ with this property exists for any family $\pi:\F\to S$
only if $g$ is even (see \cite{MR} and \cite[Theorem 3.5]{GV}).

Finally, let us provide a functorial description for a generator of the Picard group
$\Pic(\D)$.

\begin{proposition}\label{D-pic-explicit}
The image of the line bundle $\O(1)$ under the natural map \\
$\Pic(\B)\to\Pic(\D)$
generates the Picard group $\Pic(\D)$ and associates
with an object $(\C\stackrel{p}{\to}S,D)$ from $\D(S)$ the line bundle\\
$p_*$ $\left(\omega_{\C/S}^{-(g+1)}(-D)\right)$ $\in \Pic(S)$.
\end{proposition}

\begin{proof}
We compute the image under the natural map $\Pic(\D)\to\Pic(\H)$ of the element
described in the proposition. Suppose that a family $p:\C\to S$ corresponds to a
family $\pi:\F\to S$ of hyperelliptic curves. It follows
from formulas (\ref{f2}) and (\ref{f4}) that
$$
f^*(\omega_{\mathcal{C}/S}^{g+1}(D))
\cong\omega_{\mathcal{F}/S}^{g+1}(-(g+1)W)\otimes\O_{\mathcal{F}}(2W)
=\omega_{\mathcal{F}/S}^{g+1}(-(g-1)W).
$$
Combining Theorem \ref{D-pic}, Remark
\ref{comparison}, Theorem \ref{gen} and the fact that the line bundle
$\omega_{\mathcal{C}/S}^{g+1}(D)$ is trivial on any geometric fiber of $p$,
we get the desired statement.
\end{proof}

\begin{remark}
It is possible to prove Proposition \ref{D-pic-explicit} for $H_g^0$
instead of $\D$ directly, without using the stack description,
Theorem \ref{D-pic}, and Theorem \ref{gen},
see \cite[Theorem 6.3 and Remark 6.8]{GV}.
\end{remark}

\end{document}